\documentclass[11pt]{amsart}
\usepackage{amsmath} 
\usepackage{latexsym}
\usepackage{geometry}                % See geometry.pdf to learn the layout options. There are lots.
\usepackage[parfill]{parskip}    % Activate to begin paragraphs with an empty line rather than an indent
\usepackage{graphicx}
\usepackage{amssymb}
\usepackage{amsmath}
\usepackage{epstopdf}
\usepackage{hyperref}
\hypersetup{colorlinks,% 
citecolor=black,% 
filecolor=black,% 
linkcolor=black,% 
urlcolor=black,% 
pdftex}
\newtheorem{thm}{Theorem}[section]
\newtheorem{cor}[thm]{Corollary}
\newtheorem{lem}[thm]{Lemma}

\title{Division Polynomials For Twisted Edwards Curves}
\author{Richard Moloney 
and
Laura Hitt
and
Gary McGuire}
\address{
School of Mathematical Sciences\\
University College Dublin\\
Ireland\\ 
}
\email{richard.moloney@ucd.ie}
\thanks{Research supported by Claude Shannon Institute,
Science Foundation Ireland Grant 06/MI/006, and Grant 07/RFP/MATF846, 
and the Irish Research Council
for Science, Engineering and Technology}

\begin{document}

\maketitle

\begin{abstract}
This paper presents division polynomials for twisted Edwards curves. Their chief property is that they characterise the $n$-torsion points of a given twisted Edwards curve. We also present results concerning the coefficients of these polynomials, which may aid computation.
\end{abstract}

\section{Introduction}

%\subsection{Twisted Edwards curves}

Edwards  \cite{Edwards} introduced an addition law on the curves $x^2+y^2 = c^2(1+x^2y^2)$ 
for $c \in k$, where $k$ is a  field of characteristic not equal to 2. He showed that every elliptic curve over $k$ is birationally equivalent (over some extension of $k$) to a curve of this form.

In \cite{BL}, Bernstein and Lange generalised this addition law to the curves $x^2+y^2 = 1+ dx^2y^2$ for $d \in k \setminus \{0,1\}$. 
More generally, they consider $x^2 + y^2 = c^2(1+dx^2y^2)$, however, any such curve is isomorphic to one of the form  $x^2+y^2 = 1+ d^{\prime}x^2y^2$ for some $d^{\prime} \in k$, so we will assume $c = 1$. These curves are referred to as Edwards curves. Bernstein and Lange showed that if $k$ is finite, a large class of elliptic curves over $k$ (all those which have a point of order 4) can be represented in Edwards form.

In \cite{Twisted}, Bernstein et al. introduced the twisted Edwards curves $ax^2+y^2 = 1+dx^2y^2$ (where $a, \ d \in k$ are distinct and non-zero) and showed that every elliptic curve with a representation in Montgomery form is birationally equivalent to a twisted Edwards curve.

%\subsection{Aims of this paper} 
In this paper we describe a sequence of rational functions, and consequently a sequence of
polynomials, defined on the function field of a twisted Edwards curve which are analogous to the division polynomials for elliptic curves in Weierstrass form. 
 In particular, these 
polynomials characterise the $n$-torsion points of the 
twisted Edwards curve for  a positive integer $n$ (see Corollary \ref{divformcor}
and Corollary \ref{coronevar}).
These twisted Edwards division polynomials are polynomials in $y$ with coefficients in $\mathbb{Z}[a,d]$, and have degree in $y$ less than $n^2/2$.

In Theorem \ref{uniq1}
we prove a uniqueness form for elements of the function field of an Edwards curve,
analagous to the known result that elements of the function field of a
Weierstrass curve can be written uniquely in the form $p(x)+yq(x)$.
Our division polynomials (actually rational functions) are presented in this unique form.

Furthermore, we show in Section \ref{further}
 that the coefficients of a given twisted Edwards division polynomial exhibit a certain symmetry, which may reduce the amount of computation necessary for finding that polynomial.

\section{Division polynomials for Weierstrass Curves} 

We recall the division polynomials for Weierstrass curves here.

First we recall the definition of the function field of an (affine) algebraic variety. 
If $V/k$ is a variety in affine $n$-space, 
$I(V)$ denotes the ideal generated by the polynomials in $k[x_1,\dots, x_n]$
that vanish on $V$.
The affine coordinate ring of $V$ is the integral domain \[k[V]:={{k[x_1,\dots, x_n]}\over{I(V)}}.\] The function field of $V$ over $k$, denoted by $k(V)$, is defined to be the quotient field of $k[V]$.

For example, if $W$ is an elliptic curve with Weierstrass equation $v^2 = u^3+Au+B$, the function field of $W$, $k(W)$, is the quotient field of $k[u, v]/(v^2-u^3-Au-B)$.

We use $(u,v)$ as the coordinates for a curve in Weierstrass form and reserve $(x,y)$ for (twisted) Edwards curves.

If $char(k) \neq 2$ or 3, given an elliptic curve over $k$ in short Weierstrass form 
\[
W: v^2 = u^3+Au+B
\] 
with identity $\mathcal{O}$ , the division polynomials $\Psi_{n}$ are polynomials defined on the function field of $W$ for each $n \in \mathbb{N}$ by the following recursion:
\begin{align*}\Psi_{0}(u, v) &= 0\\
\Psi_{1}(u, v) &= 1\\
\Psi_{2}(u, v) &= 2v\\
\Psi_{3}(u, v) &= 3u^4+6Au^2+12Bu - A^2\\
\Psi_{4}(u, v) &= 4v(u^6+5Au^4+20Bu^3-5A^2u^2-4ABu-A^3-8B^2)\\
\Psi_{2m+1}(u, v) &= \Psi_{m+2}(u,v)\Psi_{m}^3(u,v) - \Psi_{m-1}(u,v)\Psi_{m+1}^3(u,v) \text{ for } m \geq 2\\
\Psi_{2m}(u,v) &= {\Psi_{m}(u,v)\over\Psi_{2}(u,v)}\left(\Psi_{m+2}(u,v)\Psi_{m-1}^2(u,v) -\Psi_{m-2}(u,v)\Psi_{m+1}^2(u,v)\right)\quad \text{ for }m \geq 3.\end{align*}

The $\Psi_{n}$ are polynomials in $u$ and $v$ with coefficients in $\mathbb{Z}[A,B]$. The principal properties of the division polynomials are that $\Psi_{n}(u,v) = 0$ precisely when $(u,v)$ is an $n$-torsion point of $W$ (i.e. $[n](u,v) = \mathcal{O}$), and that the multiplication-by-$n$ map $[n]:W \rightarrow W$ is characterised by the division polynomials as 
\[ [n](u,v) = \left({{u\Psi_{n}^2(u,v) - \Psi_{n-1}(u, v)\Psi_{n+1}(u, v)}\over{\Psi_{n}^2(u, v)}}, {{\Psi_{2n}(u, v)}\over{2\Psi_{n}^4(u, v)}}\right)\] 
(see e.g. \cite{Wash}, Chapters 3 , 9, \cite{Silv}, Chapter 3). 
If $n$ is odd then  $\Psi_{n} \in \mathbb{Z}[u,A,B]$, and $\Psi_n$
has degree  $(n^2-1)/2$ in $u$.  If $n$ is even then
$\Psi_{n} \in v\mathbb{Z}[u,A,B]$ with degree  $(n^2-4)/2$ in $u$.
We prove analagous results for twisted Edwards curves.

\section{Twisted Edwards Curves}

Let $k$ be a field with characteristic $\neq 2$ or $3$. 
Let $K$ be an extension field of $k$.
Let $E(K)$ be the twisted Edwards curve over $K$
with coefficients $a$ and $d$, where $a$ and $d$ are distinct and non-zero:  
\[E(K): ax^2+y^2 = 1+dx^2y^2.
\]
Points on $E(K)$ may be added by the rule 
\[(x_1, y_1)+(x_2, y_2) = \left({{x_1y_2+x_2y_1}\over{1+dx_1x_2y_1y_2}} , {{y_1y_2-ax_1x_2}\over{1-dx_1x_2y_1y_2}}\right)\]
and under this operation, the points on $E(K)$ 
form an abelian group.
The identity is $(0,1)$, and the additive inverse of a point $(x,y)$ is $(-x,y)$.
The projective closure of $E$ has singularities at  $(1:0:0)$ and $(0:1:0)$.

The twisted Edwards curve
$E(K)$ is birationally equivalent to the Weierstrass-form elliptic curve
\[W(K):  v^2 = u^3 - {{(a^2+14ad+d^2)}\over{48}}u - {{(a^3-33a^2d-33ad^2+d^3)}\over{864}}\] under the transformation 
\[
u := {{(5a-d)+(a-5d)y}\over{12(1-y)}} \ , \  v := {{(a-d)(1+y)}\over{4x(1-y)}} \quad \text{if } x(1-y)\neq 0, \]
otherwise
\begin{align*}(x,y) = (0,1) &\Rightarrow (u, v) = \mathcal{O}\\
(x,y) = (0,-1) &\Rightarrow (u, v) = \left({{a+d}\over{6}}, 0\right).\end{align*}

The inverse transformation is given by
\[x = {{6u-(a+d)}\over{6v}}, \  y = {{12u+d-5a}\over{12u+a-5d}} \quad \text{if } v(12u+a-5d)\neq 0  \]
and
\begin{align*}
(u,v) = \mathcal{O} &\Rightarrow (x,y) = (0,1)\\
(u,v) = \left({{a+d}\over{6}}, 0\right) &\Rightarrow (x,y) = (0,-1).
\end{align*}

There are 4 points on $W(\overline{k})$ that are not mapped to any point on the twisted Edwards curve. These are $(u, v) = \left({{5d-a}\over{12}} , \pm{{s(d-a)}\over{4}}\right)$ and $(u, v) =\left({{-(a+d) \pm 6t}\over{12}},0\right) $ where $s, t \in \bar{k}$ such that $s^2 = d, t^2 = ad$. We note that $\left({{-(a+d) \pm 6t}\over{12}},0\right)$ are points of order 2 on $W$, and $\left({{5d-a}\over{12}} , \pm{{s(d-a)}\over{4}}\right)$ are points of order 4 on $W$. Had we defined the birational equivalence between the projective closures of $W$ and $E$,  the points $\left(5d-a : \pm3s(d-a): 12\right)$ would map to the singular point (0:1:0) of $E$, while the points $\left(-(a+d) \pm 6t: 0: 12\right)$ would map to the singular point (1:0:0) of $E$.

\section{Function Field of a Twisted Edwards Curve}

For Weierstrass curves 
$v^2=u^3+Au+B$
it is well known (see \cite{Silv} for example) that an
element of the function field can be written uniquely in the form
\[
p(u)+vq(u)
\]
where $p(u) ,q(u)$ are polynomials in $u$.

We prove an analogous result for twisted Edwards curves $E$.
Not surprisingly, rational functions are needed in place of the polynomials.
We use the notation $\text{ord}_P (f)$ to denote the valuation of a function
$f\in K(E)$ at a point $P$.

\begin{thm} \label{uniq1}
Any function $g \in K(E)$ can be written uniquely as \[ g(x, y) = p(y) + xq(y)\] where 
$p(y)$, $q(y)$ are rational functions in $y$.
\end{thm}

\emph{Proof:} 
Let $f(x,y)=0$ be the equation defining $E$, where \[ f(x,y) = ax^2+y^2-1-dx^2y^2. \]
In $K(E)$ we have
\[
x^2=\frac{1-y^2}{a-dy^2}.
\]
It is then clear that $g$ can be written in the stated form.

Suppose this expression for $g$ is not unique. 
Then $A(y) + xB(y) = 0$ for some nonzero rational functions $A(y)$, $B(y)$. 
So \[ x = -{{A(y)}\over{B(y)}}\]  which implies
\begin{equation}\label{ords}
 \text{ord}_{(0,1)}x = \text{ord}_{(0,1)}A(y)-\text{ord}_{(0,1)}B(y).
 \end{equation}

 We obtain our contradiction by showing that the right-hand side of 
equation (\ref{ords}) is even,
but the left-hand side is equal to 1.

First we need ord$_{(0,1)}(y-1)$, so we compute 
\begin{align*}
f(x, y+1) &= ax^2+(y+1)^2-1-dx^2(y+1)^2\\
&=ax^2+y^2+2y-dx^2y^2-2dx^2y-dx^2.
\end{align*}
Then
\[ f(x, 0+1) = (a-d)x^2 \] 
which implies ord$_{(0,1)}(y-1) = 2$
since $x$ is a local uniformizer.

To find ord$_{(0,1)} x$, we use the fact that (again translating to the origin)
 \[
x^2 (a-d(y+1)^2)=y(-y-2).
\]
Note that $\text{ord}_{(0,0)} (a-d(y+1)^2)=0$ because 
 $y^2=a/d$ in the (usual) curve equation implies $a=d$, which is not allowed.
 Thus (after translation)
 \[
 \text{ord}_{(0,0)} (x^2)=  \text{ord}_{(0,0)} y
 \]
 which implies (before translation)
$\text{ord}_{(0,1)}x = 1$.

When computing ord$_{(0,1)}A(y)$, we write 
$A(y+1) = {{a(y)}\over{b(y)}}$ for some polynomials $a(y), \ b(y)$. 
Then \[\text{ord}_{(0,1)}A(y)= \text{ord}_{(0,0)}a(y)-\text{ord}_{(0,0)} b(y).\]

Let 
%\[ \frac{a(y)}{b(y)}
%= {{y^n+a_{n-1}y^{n-1}+ \dots + a_{1}y+a_{0}}\over{y^m+b_{m-1}y^{m-1}+ \dots + b_{1}y+b_{0}}} \]
%we let
$n_{0}$ be the degree of the term of smallest degree in $a(y)$,  
and similarly let $m_{0}$ be the degree of the term of smallest degree in $b(y)$. 
Then
$\text{ord}_{(0,0)}a(y) = \left(\text{ord}_{(0,0)}y\right) n_{0} = 2n_{0} $, and
similarly, $\text{ord}_{(0,0)}b(y) =  2m_{0}$. Thus $\text{ord}_{(0,1)}A(y)= 2(n_{0}-m_{0})$,
which is even.

Similarly,    $\text{ord}_{(0,1)}B(y)$ is even.
This proves that the right-hand side of (\ref{ords}) is even,
and we are done.
%Now we can say that
%\begin{align*}
%-\text{ord}_{(0,1)}x &= \text{ord}_{(0,1)}A(y) -\text{ord}_{(0,1)}B(y) \\ \Rightarrow \ -1 &=2(n_{0}-m_{0})-2(r_{0}-s_{0}).
%\end{align*}
%The left hand side is odd whie the right hand side is even,which is a contradiction.
%Thus $g(x, y)$ can be written as $g(x, y) = u(y) + xv(y)$ uniquely.
\hfill $\Box$
\bigskip

\begin{cor} \label{uniq2}
Any function $g \in K(E)$ can be written uniquely as 
\[ g(x, y) = p^{\prime}(y) + {{1}\over{x}}q^{\prime}(y)\] where 
$p^{\prime}(y)$, $q^{\prime}(y)$ are rational functions in $y$.
\end{cor}

\emph{Proof:} 
This follows from the Theorem \ref{uniq1}, and
the fact that \[x={{1}\over{x}}\cdot {{1-y^2}\over{a-dy^2}}\] on the function field of $E$.
In fact $p^\prime (y)$ is equal to $p(y)$, using the notation of Theorem \ref{uniq1},
and 
\[
q^\prime (y) = {{1-y^2}\over{a-dy^2}} \ q(y).
\]
\hfill $\Box$

\section{Division Polynomials on Twisted Edwards Curves}

We define the following rational functions $\psi_n(x,y)$ on the function field of $E$ 
recursively for $n\geq0$:

\begin{align*}\psi_0(x,y) &:= 0\\
\psi_1(x,y) &:= 1\\
\psi_2(x,y) &:= {{(a-d)(1+y)}\over{x(2(1-y))}}\\
\psi_3(x,y) &:= {{(a-d)^3(a+2ay-2dy^3-dy^4)}\over{(2(1-y))^4}}\\
\psi_4(x,y) &:= {{2(a-d)^6y(1+y)(a-dy^4)}\over{x((2(1-y))^7}}\\
\psi_{2m+1}(x,y) &:= \psi_{m+2}(x,y)\psi_m^3(x,y) - \psi_{m-1}(x,y)\psi_{m+1}^3(x,y) \quad \text{for } m \geq 2\\
\psi_{2m}(x,y) &:= {\psi_m(x,y)\over\psi_2(x,y)}\left(\psi_{m+2}(x,y)\psi_{m-1}^2(x,y) - \psi_{m-2}(x,y)\psi_{m+1}^2(x,y)\right) \quad \text{for } m \geq 3.\end{align*}

These functions are not defined at the points $(0,1)$ and $(0,-1)$.
We point out that these elements of the function field $K(E)$ are 
in the unique form given in Corollary \ref{uniq2}.

For $n \geq 1$, we also define

\begin{align*}\phi_n(x,y) &: = {{(1+y)\psi_n^2(x,y)}\over{(1-y)}} - {{4\psi_{n-1}(x,y)\psi_{n+1}(x,y)}\over{(a-d)}}\\
\text{and} \quad\omega_n(x,y) &:= {{2\psi_{2n}(x,y)}\over{(a-d)\psi_n(x,y)}}.\end{align*}

Next we show that these rational functions arise in the
multiplication-by-$n$ map.

%\textbf{Theorem 2.1:} 
\begin{thm}\label{divform}
Let $(x,y)$ be a point in $E(\overline{k})\setminus \{ (0,1), (0,-1)\} $ and $n \ge 1$ an integer. Then \[[n](x,y) = \left({{\phi_n(x,y)\psi_n(x,y)}\over{\omega_n(x,y)}}, {{\phi_n(x,y) - \psi_n^2(x,y)}\over{\phi_n(x,y) + \psi_n^2(x,y)}}\right).\]
\end{thm}

\emph{Proof:} 
Compute the division polynomials for the given Weierstrass  elliptic curve,
$W: v^2 = u^3 + Au + B, \ $ where

\[A = - {{(a^2+14ad+d^2)}\over{48}},\qquad B = - {{(a^3-33a^2d-33ad^2+d^3)}\over{864}}.\]
We get
\begin{align*}\Psi_{0}(u, v) &= 0\\
\Psi_{1}(u, v) &= 1\\
\Psi_{2}(u, v) &= 2v\\
\Psi_{3}(u, v) &= 3u^4+6Au^2+12Bu - A^2\\
\Psi_{4}(u, v) &= 4v(u^6+5Au^4+20Bu^3-5A^2u^2-4ABu-A^3-8B^2)\\
\Psi_{2m+1}(u, v) &= \Psi_{m+2}(u,v)\Psi_{m}^3(u,v) - \Psi_{m-1}(u,v)\Psi_{m+1}^3(u,v) \text{ for } m \geq 2\\
\Psi_{2m}(u,v) &= {\Psi_{m}(u,v)\over\Psi_{2}(u,v)}\left(\Psi_{m+2}(u,v)\Psi_{m-1}^2(u,v) - \Psi_{m-2}(u,v)\Psi_{m+1}^2(u,v)\right)\quad \text{ for }m \geq 3.\end{align*}

Substituting \begin{align*}&A = - {{(a^2+14ad+d^2)}\over{48}},\quad B = - {{(a^3-33a^2d-33ad^2+d^3)}\over{864}}\quad \text{and}\\
&u := {{(5a-d)+(a-5d)y}\over{12(1-y)}},   \qquad v := {{(a-d)(1+y)}\over{4x(1-y)}},\end{align*}

for the cases $0,1,2,3,4$ we see that $\Psi_{i}(u, v) = \psi_i(x,y)$ for $i = 0,1,2,3,4$. Hence, as the recursion relations for the two sets of functions $\Psi_{i}(u, v)$ and $\psi_i(x,y)$ 
are identical for $i \geq 5$, we have that $\Psi_{n}(u, v) = \psi_n(x,y)$ for all integers $n \geq 0$.

From here on we will use the abbreviated notations $\psi_n$ for $\psi_n(x,y)$, $\phi_n$ for $\phi_n(x,y)$ and $\omega	_n$ for $\omega_n(x,y)$.
Let $\left(x_n , y_n\right) = [n](x,y)$, and $\left(u_n, v_n\right) = [n]_{W}\left(u, v\right)$.

From the properties of the division polynomials,

\[u_n = u -{{\Psi_{n-1}(u, v)\Psi_{n+1}(u, v)}\over{\Psi_{n}^2(u, v)}}, \quad
v_n = {{\Psi_{2n}(u, v)}\over{2\Psi_{n}^4(u, v)}},\]

i.e.,

\[u_n = u -{{\psi_{n-1}\psi_{n+1}}\over{\psi_{n}^2}}, \quad
v_n = {{\psi_{2n}}\over{2\psi_{n}^4}},\]

and, applying the birational equivalence gives

\[x_n = {{6u_n - (a+d)}\over{6v_n}}, \quad y_n =  {{12u_n+d-5a}\over{12u_n+a-5d}},\]

\begin{align*}x_n &= {{2\psi_n^4}\over{\psi_{2n}}}\left({{5a-d+(a-5d)y}\over{12(1-y)}} - {{\psi_{n-1}\psi_{n+1}}\over{\psi_n^2}}- {{a+d}\over{6}}\right)\\
&= {{\psi_n^2}\over{\psi_{2n}}}\left({{(a-d)(1+y)\psi_n^2}\over{2(1-y)}}- {2\psi_{n-1}\psi_{n+1}}\right)\end{align*}

while \begin{align*} {{\phi_n\psi_n}\over{\omega_n}} &= {{(a-d)\psi_n^2}\over{2\psi_{2n}}}\left(\left({{1+y}\over{1-y}}\right)\psi_n^2-{{4\psi_{n-1}\psi_{n+1}}\over{a-d}}\right)\\
&= {{\psi_n^2}\over{\psi_{2n}}}\left({{(a-d)(1+y)\psi_n^2}\over{2(1-y)}}- {2\psi_{n-1}\psi_{n+1}}\right)\\
 &=x_n.\\
\end{align*}
Also,
\[y_n = {{12u_n+d-5a}\over{12u_n+a-5d}}\]
and
\begin{align*}12u_n+d-5a &= {{5a-d+(a-5d)y}\over{(1-y)}} - 12{{\psi_{n-1}\psi_{n+1}}\over{\psi_n^2}} +d-5a\\
&={{6(a-d)y}\over{1-y}}- 12{{\psi_{n-1}\psi_{n+1}}\over{\psi_n^2}}\\
12u_n+a-5d &= {{6(a-d)}\over{1-y}}- 12{{\psi_{n-1}\psi_{n+1}}\over{\psi_n^2}}\end{align*}

so \[y_n = {{(a-d)y\psi_n^2-2(1-y)\psi_{n-1}\psi_{n+1}}\over{(a-d)\psi_n^2-2(1-y)\psi_{n-1}\psi_{n+1}}}\]

and \begin{align*}{{\phi_n-\psi_n^2}\over{\phi_n+\psi_n^2}} &= {{\left({{1+y}\over{1-y}}\right)\psi_n^2-{{4\psi_{n-1}\psi_{n+1}}\over{a-d}}-\psi_n^2}\over{\left({{1+y}\over{1-y}}\right)\psi_n^2-{{4\psi_{n-1}\psi_{n+1}}\over{a-d}}+\psi_n^2}}\\
&={{(a-d)y\psi_n^2-2(1-y)\psi_{n-1}\psi_{n+1}}\over{(a-d)\psi_n^2-2(1-y)\psi_{n-1}\psi_{n+1}}}\\
& = y_n .\end{align*}

Hence \[[n](x,y) = \left({{\phi_n(x,y)\psi_n(x,y)}\over{\omega_n(x,y)}}, {{\phi_n(x,y) - \psi_n^2(x,y)}\over{\phi_n(x,y) + \psi_n^2(x,y)}}\right).\] 
\hfill $\Box$

\begin{cor}\label{divformcor}
Let $P=(x,y)$ be in
$E(\overline{k})\setminus \{ (0,1), (0,-1)\} $ and let $n\geq 1$.
Then $P$ is an $n$-torsion point of $E$ if and 
only if $\psi_n(P) = 0$.
\end{cor}
  
\emph{Proof:}  Since the identity is $(0,1)$,
the result is clear from Theorem \ref{divform}.
\hfill $\Box$

\bigskip

 So the $\psi_n(x,y)$, though they are rational functions, can be seen as analogues of division polynomials.  Here are the first seven $\psi_n(x,y)$:

\begin{align*}
\psi_0 &= 0\\
\psi_1 &= 1\\
\psi_2 &= {{(a-d)(y+1)}\over{x(2(1-y))}}\\
\psi_3 &= {{(a-d)^3(-dy^4-2dy^3+2ay+a)}\over{(2(1-y))^4}}\\
\psi_4 &= {{2(a-d)^6(-2dy^6-2dy^5+2ay^2+2ay)}\over{x((2(1-y))^7}}\\
\psi_5 &={{(a-d)^9(d^3y^{12}+5d^3y^{11}+\dots-5a^3y-a^3)}\over{(2(1-y))^{12}}}\\
\psi_6 &= {{(a-d)^{13}(-d^4y^{17}-d^4y^{16}+(4ad^3+4d^4)y^{15}+\dots+(4a^3d+4a^4)y^2-a^4y-a^4)}\over{x((2(1-y))^{17}}}.
\end{align*}

As we said earlier,  these elements of the function field $K(E)$ are 
in the unique form given in Corollary \ref{uniq2}.

The apparent patterns here are proved in the next theorem.

\section{Division Polynomials}

The next theorem isolates the key polynomial in the numerator of $\psi_n$,
which we call $\tilde{\psi} (y)$.
These polynomials could also be called the division polynomials for twisted Edwards curves.

\begin{thm} \label{main}
We have
\begin{equation*}
\psi_n(x,y)= \left\{
	\begin{array}{rl}
		{(a-d)^{k(n)}\tilde{\psi}_{n}(y)/(2(1-y))^{m(n)}} & \text{ if $n$ is odd}\\
		\\
		{(a-d)^{k(n)}\tilde{\psi}_{n}(y)/x(2(1-y))^{m(n)}} & \text{ if $n$ is even}
	\end{array} \right.
\end{equation*}
where
\begin{equation*}
m(n)= \left\{
	\begin{array}{rl}
		{{n^2-1}\over{2}} & \text{ if n is odd}\\
		\\
		{{n^2-2}\over{2}} & \text{ if n is even}
	\end{array} \right.
\end{equation*}
and
\[
k(n) = \left\lfloor{{3n^2}\over{8}}\right\rfloor\]
and
\begin{align*}
\tilde{\psi}_0(y) &= 0\\
\tilde{\psi}_1(y) &= 1\\
\tilde{\psi}_2(y) &= y+1\\
\tilde{\psi}_3(y) &= -dy^4-2dy^3+2ay+a\\
\tilde{\psi}_4(y) &= -2y(y+1)(dy^4-a) = -2dy^6-2dy^5+2ay^2+2ay,\end{align*}
and
\begin{equation*}
\tilde{\psi}_{2r+1}(y) = \left\{
	\begin{array}{rl}
		{{4(a-d)(a-dy^2)^2\tilde{\psi}_{r+2}(y)\tilde{\psi}_{r}^3(y)}\over{(y+1)^2}} - \tilde{\psi}_{r-1}(y)\tilde{\psi}_{r+1}^3(y)& \text{ if } r\equiv 0 \pmod{4},\ r\geq4\\
		{\tilde{\psi}_{r+2}(y)\tilde{\psi}_{r}^3(y)} - {{4(a-dy^2)^2\tilde{\psi}_{r-1}(y)\tilde{\psi}_{r+1}^3(y)}\over{(y+1)^2}}& \text{ if } r\equiv 1 \pmod{4},\ r\geq5\\
		{{4(a-dy^2)^2\tilde{\psi_{r+2}}(y)\tilde{\psi_{r}}^3(y)}\over{(y+1)^2}} - \tilde{\psi}_{r-1}(y)\tilde{\psi}_{r+1}^3(y)& \text{ if } r\equiv 2 \pmod{4},\ r\geq2\\
		{\tilde{\psi}_{r+2}(y)\tilde{\psi}_{r}^3(y)} - {{4(a-d)(a-dy^2)^2\tilde{\psi}_{r-1}(y)\tilde{\psi}_{r+1}^3(y)}\over{(y+1)^2}}& \text{ if } r\equiv 3 \pmod{4},\ r\geq3
	\end{array} \right.
\end{equation*}
and
\begin{equation*}
\tilde{\psi}_{2r}(y) = \left\{
	\begin{array}{rl}
		{{\tilde{\psi}_{r}(y)}\over{y+1}}\left(\tilde{\psi}_{r+2}(y)\tilde{\psi}_{r-1}^2(y) - \tilde{\psi}_{r-2}(y)\tilde{\psi}_{r+1}^2(y)\right)& \text{ if } r\equiv 0 \pmod{4},\ r\geq4\\
		{{\tilde{\psi}_{r}(y)}\over{y+1}}\left((a-d)\tilde{\psi}_{r+2}(y)\tilde{\psi}_{r-1}^2(y) - \tilde{\psi}_{r-2}(y)\tilde{\psi}_{r+1}^2(y)\right)& \text{ if } r\equiv 1 \pmod{4},\ r\geq5\\
		{{\tilde{\psi}_{r}(y)}\over{y+1}}\left(\tilde{\psi}_{r+2}(y)\tilde{\psi}_{r-1}^2(y) - \tilde{\psi}_{r-2}(y)\tilde{\psi}_{r+1}^2(y)\right)& \text{ if } r\equiv 2 \pmod{4},\ r\geq6\\
		{{\tilde{\psi}_{r}(y)}\over{y+1}}\left(\tilde{\psi}_{r+2}(y)\tilde{\psi}_{r-1}^2(y) - (a-d)\tilde{\psi}_{r-2}(y)\tilde{\psi}_{r+1}^2(y)\right)& \text{ if } r\equiv 3 \pmod{4},\ r\geq3.
		\end{array} \right.
\end{equation*}
\end{thm}

\emph{Proof:} 

First observe for all $ t \in \mathbb{Z}, \ t>0$,
\begin{align*}m(4t) &= {{16t^2-2}\over{2}} = 8t^2-1\\
m(4t\pm1)&={{(4t\pm1)^2-1}\over{2}}={{16t^2\pm8t}\over{2}}=8t^2\pm4t\\
m(4t\pm2)&={{(4t\pm2)^2-2}\over{2}}={{16t^2\pm16t+2}\over{2}} = 8t^2\pm8t+1\\
m(4t\pm3)&={{(4t\pm3)^2-1}\over{2}}={{16t^2\pm24t+8}\over{2}} = 8t^2\pm12t+4
\end{align*}
and
\begin{align*}
k(4t) &= \left\lfloor{{3(4t)^2}\over{8}}\right\rfloor = \left\lfloor{6t^2}\right\rfloor = 6t^2\\
k(4t\pm1) &= \left\lfloor{{{3(4t\pm1)^2}\over{8}}}\right\rfloor =\left\lfloor{6t^2\pm3t+{{3}\over{8}}}\right\rfloor= 6t^2\pm3t\\
k(4t\pm2) &= \left\lfloor{{{3(4t\pm2)^2}\over{8}}}\right\rfloor =\left\lfloor{6t^2\pm6t+{{12}\over{8}}}\right\rfloor= 6t^2\pm6t+1\\
k(4t\pm3) &= \left\lfloor{{{3(4t\pm3)^2}\over{8}}}\right\rfloor =\left\lfloor{6t^2\pm9t+{{27}\over{8}}}\right\rfloor= 6t^2\pm9t+3.
\end{align*}

The proof is by induction.
The claim is true for $n = 0 \dots 4$.

Assume true for $0 \dots n-1$

\underline{Case 1:} $n \equiv 0 \pmod{8}$ i.e. $ n = 8l$ for some $l \in \mathbb{Z}$. Let $r = 4l$.

By definition,
\begin{align*}\psi_{n} &= {{\psi_r}\over{\psi_2}}\left(\psi_{r+2}\psi_{r-1}^2 - \psi_{r-2}\psi_{r+1}^2\right)\\
&= {{(a-d)^{k(r)-1}\tilde{\psi}_{r}}\over{(y+1)(2(1-y))^{m(r)-1}}}\left({{(a-d)^{k(r+2)+2k(r-1)}\tilde{\psi}_{r+2}\tilde{\psi}_{r-1}^2}\over{x(2(1-y))^{m(r+2)+2m(r-1)}}} - {{(a-d)^{k(r-2)+2k(r+1)}\tilde{\psi}_{r-2}\tilde{\psi}_{r+1}^2}\over{x(2(1-y))^{m(r-2)+2m(r+1)}}}\right).\end{align*}

Also,
\begin{align*}m(4l)-1+m(4l+2) + 2m(4l-1) &= 8l^2-1-1+8l^2+8l+1+16l^2-8l\\
&=32l^2-1 = m(8l) = m(n)\\
m(4l)-1+m(4l-2) + 2m(4l+1) &= 8l^2-1-1+8l^2-8l+1+16l^2+8l\\
&=32l^2-1 = m(8l) = m(n)\end{align*}

and
\begin{align*}k(4l)-1+k(4l+2) + 2k(4l-1) &= 6l^2-1+6l^2+6l+1+12l^2-6l\\
&=24l^2 = k(8l) = k(n)\\
k(4l)-1+k(4l-2) + 2k(4l+1) &= 6l^2-1+6l^2-6l+1+12l^2+6l\\
&=24l^2 = k(8l) = k(n).
\end{align*}

So \begin{align*}
\psi_n &={{(a-d)^{k(n)}}\over{x(y+1)(2(1-y))^{m(n)}}}\left(\tilde{\psi}_{r}\left(\tilde{\psi}_{r+2}\tilde{\psi}_{r-1}^2 - \tilde{\psi}_{r-2}\tilde{\psi}_{r+1}^2\right)\right)\\
&={{(a-d)^{k(n)}\tilde{\psi}_{n}(y)}\over{x(2(1-y))^{m(n)}}} \ .
\end{align*}

\underline{Case 2:}  $n \equiv 1 \pmod{8}$ i.e. $ n = 8l+1$ for some $l \in \mathbb{Z}$. Let $r = 4l$.

By definition
\begin{align*}\psi_{n} &= \psi_{r+2}\psi_r^3 - \psi_{r-1}\psi_{r+1}^3\\
&={{(a-d)^{k(r+2)+3k(r)}\tilde{\psi}_{r+2}\tilde{\psi}_r^3}\over{y^4(2(1-x))^{m(r+2)+3m(r)}}} - {{(a-d)^{k(r-1)+3k(r+1)}\tilde{\psi}_{r-1}\tilde{\psi}_{r+1}^3}\over{(2(1-y))^{m(r-1)+3m(r+1)}}}.
\end{align*}

Using the curve equation \[ax^2+y^2  = 1 + dx^2y^2\]
gives\begin{align*}&x^2 = {{(1-y^2)}\over{(a-dy^2)}} = {{(1-y)(1+y)}\over{(a-dy^2)}}\\
\Rightarrow &x^4= {{(1-y)^2(1+y)^2}\over{(a-dy^2)^2}}\end{align*}

so \[\psi_n = {{4(a-d)^{k(r+2)+3k(r)}(a-dy^2)^2\tilde{\psi}_{r+2}\tilde{\psi}_r^3}\over{(y+1)^2(2(1-y))^{m(r+2)+3m(r)+2}}} - {{(a-d)^{k(r-1)+3k(r+1)}\tilde{\psi}_{r-1}\tilde{\psi}_{r+1}^3}\over{(2(1-y))^{m(r-1)+3m(r+1)}}}\ .\]

Again,
\begin{align*}k(4l+2)+3k(4l) & =6l^2+6l+1 + 18l^2 = 24l^2+6l+1\\& = k(n)+1\\
k(4l-1)+3k(4l+1)&= 6l^2-3l+18l^2+9l=24l^2+6l \\&= k(n)
\end{align*}

and
\begin{align*}
m(4l+2)+3m(4l)+2 &= 8l^2+8l+1+24l^2-3 +2= 32l^2 +8l\\&=m(n)\\
m(4l-1)+3m(4l+1) &=8l^2-4l+24l^2+12l = 32l^2 +8l\\&=m(n).
\end{align*}
Hence\[\psi_n = {{4(a-d)(a-dy^2)^2\tilde{\psi}_{r+2}(y)\tilde{\psi}_{r}^3(y)}\over{(y+1)^2}} - \tilde{\psi}_{r-1}(y)\tilde{\psi}_{r+1}^3(y)\ .\]

\underline{Cases 3,$\dots$8:}  $n \equiv 2, \dots 7 \pmod{8}$. Similar.
\hfill $\Box$

\bigskip

\begin{cor}\label{coronevar}  
Let $P=(x,y)$ be in
$E(\overline{k})\setminus \{ (0,1)\} $ and let $n\geq 1$.
Then 
\begin{center}
$P$ is an $n$-torsion point of $E$\quad  if and only if  \quad  $\tilde\psi_n(y) = 0$.
\end{center}
\end{cor}
  
Proof: 
The result follows from Corollary \ref{divformcor} and Theorem \ref{main}.
\hfill $\Box$

\bigskip

\section{Further Facts}\label{further}

Here are some more facts about the $\tilde{\psi}$.

\begin{thm} \label{int}
$\tilde{\psi}_n(y) \in \mathbb{Z}[a,d,y] \ \ \forall n >0$, and $(y+1)$ divides $\tilde{\psi}_n(y)$ if n is even
\end{thm}

\emph{Proof:} Proof is by induction. The statement is true for $n = 0, 1, 2, 3, 4$. Now suppose it is true for $0,1,2, \dots, n-1$:

\underline{Case 1:} $n \equiv 0 \pmod{8}$ i.e. $ n = 8l$ for some $l \in \mathbb{Z}$. Let $r = 4l$.

Then $\tilde{\psi}_n(y) = {{\tilde{\psi}_{r}(y)}\over{y+1}}\left(\tilde{\psi}_{r+2}(y)\tilde{\psi}_{r-1}^2(y) - \tilde{\psi}_{r-2}(y)\tilde{\psi}_{r+1}^2(y)\right)$

and $\tilde{\psi}_{r}(y), \ \tilde{\psi}_{r+2}(y), \ \tilde{\psi}_{r-1}(y),\ \tilde{\psi}_{r-2}(y), \ \tilde{\psi}_{r+1}(y) \in \mathbb{Z}[a,d,y] $. Also, $(y+1)$ divides $\tilde{\psi}_{r}(y), \ \tilde{\psi}_{r+2}(y)$, and $\tilde{\psi}_{r-2}(y)$ by hypothesis. Hence $\tilde{\psi}_n(y) \in \mathbb{Z}[a,d,y] $ and $(y+1)$ divides $\tilde{\psi}_{n}(y)$.

\underline{Case 2:}  $n \equiv 1 \pmod{8}$ i.e. $ n = 8l+1$ for some $l \in \mathbb{Z}$. Let $r = 4l$.

Then $\tilde{\psi}_n(y) = {{4(a-d)(a-dy^2)^2\tilde{\psi}_{r+2}(y)\tilde{\psi}_{r}^3(y)}\over{(y+1)^2}} - \tilde{\psi}_{r-1}(y)\tilde{\psi}_{r+1}^3(y)$

and $\tilde{\psi}_{r+2}(y), \ \tilde{\psi}_{r}(y), \ \tilde{\psi}_{r-1}(y), \ \tilde{\psi}_{r+1}(y) \in \mathbb{Z}[a,d,y] $. Also, $(y+1)$ divides $\tilde{\psi}_{r}(y)$ and $\tilde{\psi}_{r+2}(y)$ by hypothesis. Hence $\tilde{\psi}_n(y) \in \mathbb{Z}[a,d,y]$.

\underline{Cases 3,$\dots$8:}  $n \equiv 2, \dots 7 \pmod{8}$. Similar.
\hfill $\Box$
\bigskip

Theorem \ref{leading} and Corollary \ref{degree} provide results for the degrees of these polynomials $\tilde{\psi}_n(y)$, and  Theorem \ref{symm} shows that the coefficients of the polynomials exhibit a large amount of symmetry.

\begin{thm} \label{leading} If $char(k) = 0$ or $4 \cdot char(k) \nmid n$, then $\tilde{\psi}_{n}(y)$ has leading term (term of largest degree in $y$)
\begin{equation*} \left\{
	\begin{array}{rl}
		\delta(n) d^{m(n)-k(n)}y^{m(n)}& \text{ if } n \not\equiv 0 \pmod{4}\\
		\\
		\delta(n) d^{m(n)-k(n)}y^{m(n)-1}& \text{ if } n \equiv 0 \pmod{4}
	\end{array} \right.
\end{equation*}
where
\begin{equation*} \delta (n) = \left\{
	\begin{array}{rl}
		{{n}\over{2}} & \text{ if } n \equiv 0 \pmod{8}\\
		\\
		- {{n}\over{2}} & \text{ if } n \equiv 4 \pmod{8}\\
		\\
		1 & \text{ if } n \equiv 1,2, \text{ or } 5 \pmod{8}\\
		\\
		-1& \text{ if } n \equiv 3,6, \text{ or } 7 \pmod{8}
	\end{array} \right.
\end{equation*}
and $m(n)$, $k(n)$ are as defined in Theorem \ref{main}.

If $char(k) \neq 0$ and $4\cdot char(k) \mid n$, then $deg(\tilde{\psi}_{n}(y)) < m(n)-1$ .
\end{thm}

\emph{Proof:} Proof is by induction. The statement is true for  $n = 0, 1, 2, 3, 4$. Now suppose it is true for $0,1,2, \dots, n-1$:

\underline{Case 1:} $n \equiv 0 \pmod{8}$ i.e. $ n = 8l$ for some $l \in \mathbb{Z}$. Let $r = 4l$.
Then
\begin{align*}
	\tilde{\psi}_n(y) =& {{\tilde{\psi}_{r}(y)}\over{y+1}}\left(\tilde{\psi}_{r+2}(y)\tilde{\psi}_{r-1}^2(y) - \tilde{\psi}_{r-2}(y)\tilde{\psi}_{r+1}^2(y)\right)\\
	=&(\delta(r) d^{m(r)-k(r)}y^{m(r)-2}+\dots)\times \\
	&[(\delta(r+2)(\delta(r-1))^2d^{m(r+2)+2m(r-1)-k(r+2)-2k(r-1)}y^{m(r+2)+2m(r-1)} +\dots)\\& \  -(\delta(r-2)(\delta(r+1))^2d^{m(r-2)+2m(r+1)-k(r-2)-2k(r+1)}y^{m(r-2)+2m(r+1)}+ \dots)]
\end{align*}

So, computing the $m$'s and $k$'s as in previous proofs, and noting that \begin{align*}\delta(r) = \pm2l,\ &\delta(r+2) = \pm1, \ \delta(r-1) = -1, \\ &\delta(r-2) = \mp1, \ \delta(r+1) = 1, \end{align*}
the leading term is thus 
\[ \pm2ld^{m(n)-k(n)}y^{m(r)-2}(\pm y^{m(r+2)+2m(r-1)} \pm y^{m(r-2)+2m(r+1)})\] \[ = {{n}\over{2}}d^{m(n)-k(n)}y^{m(n)-1} \]
\[=\delta(n)d^{m(n)-k(n)}y^{m(n)-1}. \]

The only exception being if $char(k) \neq 0$ and $char(k) \mid r$, (i.e. if $char(k) \mid n$) in which case,  $deg(\tilde{\psi}_{r}(y)) < m(r)-1$ and  $deg(\tilde{\psi}_{n}(y)) < m(n)-1$.

\underline{Case 2:}  $n \equiv 1 \pmod{8}$ i.e. $ n = 8l+1$ for some $l \in \mathbb{Z}$. Let $r = 4l$.

Then $\tilde{\psi}_n(y) = {{4(a-d)(a-dy^2)^2\tilde{\psi}_{r+2}(y)\tilde{\psi}_{r}^3(y)}\over{(y+1)^2}} - \tilde{\psi}_{r-1}(y)\tilde{\psi}_{r+1}^3(y)$.

The degree (in $y$) of the first term above is $m(r+2)+3(m(r)-1)+4-2 = 32l^2+8l-3$.

The degree (in $y$) of the second term is $m(r-1)+3m(r+1) = 32l^2+8l$
Thus ${{4(a-d)(a-dy^2)^2\tilde{\psi}_{r+2}(y)\tilde{\psi}_{r}^3(y)}\over{(y+1)^2}}$ does not contribute to the leading term which is 
\[-\delta(r-1)(\delta(r+1))^3d^{m(r-1)+3m(r+1)-k(r-1)-3k(r+1)}y^{32l^2+8l}. \] Now,
\[\delta(r-1) = -1, \  \delta(r+1) = 1, \ \delta(n) = 1\]
\[k(r-1)+3k(r+1) =24l^2+6l\]
\[m(n) = m(8l+1)=32l^2+8l-(24l^2+6l) = 8l^2+2l. \]
So the leading term is $d^{m(n)-k(n)}y^{m(n)} = \delta(n)d^{m(n)-k(n)}y^{m(n)}$, as required.

The only exceptional case is if $char(k) \neq 0$ and $char(k) \mid r$, in which case $deg(\tilde{\psi}_{r}(y)) < m(r)-1$, but as  $\tilde{\psi}_{r}(y)$ does not contribute to the leading term, this does not affect the result.

\underline{Cases 3,$\dots$8:}  $n \equiv 2, \dots 7 \pmod{8}$. Similar.
\hfill $\Box$
\bigskip

\begin{cor} \label{degree}  If $4 \nmid n$, then $deg(\tilde{\psi}_{n}(y)) = m(n)$ where
\begin{equation*}
m(n)= \left\{
	\begin{array}{rl}
		{{n^2-1}\over{2}} & \text{ if n is odd}\\
		\\
		{{n^2-2}\over{2}} & \text{ if n is even.}
	\end{array} \right.
\end{equation*}
If $4 \mid n$ and $char(k) \nmid n$,  $deg(\tilde{\psi}_{n}(y)) = m(n)-1.$

Otherwise $deg(\tilde{\psi}_{n}(y))< m(n)-1.$
\end{cor}

\emph{Proof:} Immediate from Theorem \ref{leading} .   \hfill $\ \Box$

The only case where the degree of the polynomial $\tilde{\psi}_{n}$ is not known precisely is when $4 \cdot char(k) \mid n$.

\begin{lem} \label{last} If $char(k) = 0$ or $4 \cdot char(k) \nmid n$, then $\tilde{\psi}_{n}(y)$ has final term (term of least degree in $y$)
\begin{equation*} \left\{
	\begin{array}{rl}
		\epsilon(n) a^{m(n)-k(n)}& \text{ if } n \not\equiv 0 \pmod{4}\\
		\\
		\epsilon(n) a^{m(n)-k(n)}y& \text{ if } n \equiv 0 \pmod{4}
	\end{array} \right.
\end{equation*}
where
\begin{equation*} \epsilon (n) = \left\{
	\begin{array}{rl}
		-{{n}\over{2}} & \text{ if } n \equiv 0 \pmod{8}\\
		\\
		 {{n}\over{2}} & \text{ if } n \equiv 4 \pmod{8}\\
		\\
		1 & \text{ if } n \equiv 1,2, \text{ or } 3 \pmod{8}\\
		\\
		-1& \text{ if } n \equiv 5,6, \text{ or } 7 \pmod{8}
	\end{array} \right.
\end{equation*}
and $m(n)$, $k(n)$ are as defined in Theorem \ref{main}.

If $char(k) \neq 0$ and $4\cdot char(k) \mid n$, then the term of least degree has degree greater than 1. 
\end{lem}

\emph{Proof:} Similar to proof of Theorem \ref{leading}.  
\hfill $\ \Box$
\bigskip

Recall from Theorem \ref{int} that $\tilde{\psi}_{n}(y) = \tilde{\psi}_{n}(a,d,y) \in  \mathbb{Z}[a,d,y]$. If we write $\tilde{\psi}_{n}$ in the form \[ \tilde{\psi}_{n}(a,d,y) = \alpha_{m(n)}y^{m(n)} + \alpha_{m(n)-1}y^{m(n)-1}+ \dots + \alpha_{1}y+\alpha_{0}\] where $m(n)$ is as defined in Theorem \ref{main} (so, in particular, if $4 \mid n, \  \alpha_{m(n)} = \alpha_{0} =0$) and $\alpha_{i} \in \mathbb{Z}[a,d]$, then we define \[\tilde{\psi}_{n}^*(a,d,y):=\alpha_{0}y^{m(n)}+\alpha_{1}y^{m(n)-1}+\dots+\alpha_{m(n)-1}y+\alpha_{m(n)}\]

\begin{lem} \label{hom} $\tilde{\psi}_{n}(a,d,y)$, considered as a polynomial in $a$ and $d$ (with 
coefficients in $\mathbb{Z}[a,d]$) is homogeneous of degree $m(n)-k(n)$.
\end{lem}

\emph{Proof:} Proof is by induction using Theorem \ref{main}.  
\hfill $\ \Box$
\bigskip

\begin{thm} \label{symm}
Consider $\tilde{\psi}_{n}(a,d,y) \in  \mathbb{Z}[a,d,y]$, as a polynomial in three variables. Then $\tilde{\psi}_{n}(a,d,y) = \tilde{\psi}_{n}^*(-d,-a,y)$.
\end{thm}

\emph{Proof:} We can restate this theorem as: If \[\tilde{\psi}_{n}(a,d,y) = \alpha_{m(n)}(a,d)y^{m(n)} + \alpha_{m(n)-1}(a,d)y^{m(n)-1}+ \dots + \alpha_{1}(a,d)y+\alpha_{0}(a,d)\]then \[\tilde{\psi}_{n}(a,d,y) = \alpha_{0}(-d,-a)y^{m(n)}+\alpha_{1}(-d,-a)y^{m(n)-1}+\dots+\alpha_{m(n)-1}(-d,-a)y+\alpha_{m(n)}(-d,-a).\]

If $E$ is as defined at the outset, \[E: ax^2+y^2  = 1+dx^2y^2\]
and we let $E^{\prime}$ be the twisted Edwards curve \[E^{\prime}: dx^2+y^2 = 1+ax^2y^2\]
then the birational equivalence $(x,y) \mapsto \left(x, {{1}\over{y}}\right)$ maps $E$ to $E^{\prime}$, and $E^{\prime}$ to $E$.

Now, \[\psi_{n}(x,y) = {{(a-d)^{k(n)}\tilde{\psi}_{n}(y)}\over{(2(1-y))^{m(n)}x^{\gamma(n)}}}\] where
 \begin{equation*} \gamma(n) = \left\{
	\begin{array}{rl}
		1 &\text{ if $n$ is even}\\
		0 &\text{ if $n$ is odd}
	\end{array} \right.
\end{equation*}

and \[\psi_{n}^{\prime}(x,y) = {{(d-a)^{k(n)}\tilde{\psi}_{n}^{\prime}(y)}\over{(2(1-y))^{m(n)}x^{\gamma(n)}}}\]

where $\psi_{n}^{\prime}(x,y), \ \tilde{\psi}_{n}^{\prime}(y)$ are the relevant functions defined on $E^{\prime}$.

Now,
\begin{align*}
\psi_{n}^{\prime}(x, {{1}\over{y}}) &= {{(d-a)^{k(n)}\tilde{\psi}_{n}^{\prime}({{1}\over{y}})}\over{(2(1-{{1}\over{y}}))^{m(n)}x^{\gamma(n)}}}\\
&={{(a-d)^{k(n)}((-1)^{m(n)-k(n)}y^{m(n)}\tilde{\psi}_{n}^{\prime}({{1}\over{y}}))}\over{(2(1-y))^{m(n)}x^{\gamma(n)}}}
\end{align*}

and by theorem \ref{leading}, $(-1)^{m(n)-k(n)}y^{m(n)}\tilde{\psi}_{n}^{\prime}({{1}\over{y}}) \in \mathbb{Z}[a,d,y]$.
 
By the birational equivalence, for any $(x, y) \in E$, \[\psi_{n}(x,y) = 0 \Leftrightarrow \psi_{n}^{\prime}\left(x, {{1}\over{y}}\right)=0\]

so \[\tilde{\psi}_{n}(y) = 0 \Leftrightarrow (-1)^{m(n)-k(n)}y^{m(n)}\tilde{\psi}_{n}^{\prime}({{1}\over{y}}) =0\]

which gives \[\tilde{\psi}_{n}(y) = t(-1)^{m(n)-k(n)}y^{m(n)}\tilde{\psi}_{n}^{\prime}({{1}\over{y}})\]
for some $t$. By comparing leading terms using theorems \ref{leading} and \ref{last}, we get $t = 1$, i.e.,
 
\[\tilde{\psi}_{n}(y) = (-1)^{m(n)-k(n)}y^{m(n)}\tilde{\psi}_{n}^{\prime}({{1}\over{y}}).\]

Now, \[\tilde{\psi}_{n}(a,d,y) = \alpha_{m(n)}(a,d)y^{m(n)} + \alpha_{m(n)-1}(a,d)y^{m(n)-1}+ \dots + \alpha_{1}(a,d)y+\alpha_{0}(a,d)\]
and \[\tilde{\psi}^{\prime}_{n}(a,d,y) = \alpha_{m(n)}(d,a)y^{m(n)} + \alpha_{m(n)-1}(d,a)y^{m(n)-1}+ \dots + \alpha_{1}(d,a)y+\alpha_{0}(d,a).\]

Recall (lemma \ref{hom}) that each of the $\alpha_i$ is homogeneous in $a$ and $d$ of degree $m(n) -k(n)$, so 
\[(-1)^{m(n)-k(n)}\tilde{\psi}^{\prime}_{n}(a,d,y) = \alpha_{m(n)}(-d,-a)y^{m(n)} + \alpha_{m(n)-1}(-d,-a)y^{m(n)-1}+ \dots + \alpha_{1}(-d,-a)y+\alpha_{0}(-d,-a)\]

and \begin{align*}
(-1)^{m(n)-k(n)}y^{m(n)}\tilde{\psi}_{n}^{\prime}({{1}\over{y}}) &= \alpha_{m(n)}(-d,-a) + \alpha_{m(n)-1}(-d,-a)y+ \dots \\
&\quad +\alpha_{1}(-d,-a)y^{m(n)-1}+\alpha_{0}(-d,-a)y^{m(n)}\\
&=\tilde{\psi}_{n}^*(-d,-a,y).
\end{align*}
Hence, $\tilde{\psi}_{n}(a,d,y) = \tilde{\psi}_{n}^*(-d,-a,y).$ 
\hfill $\ \Box$
\bigskip

%\section{Conclusion}
%The division polynomials defined for elliptic curves in Weierstrass form map to rational functions under the transformation to a twisted Edwards curve (where it exists). In theorem \ref{main}, we showed that these rational functions are straightforwardly related to a sequence of polynomials  with integer coefficients defined by a recursion similar to (but slightly more complicated than) that which defines the division polynomials. These new polynomials have the same degree as the division polynomials if $n$ is odd. In theorem \ref{symm}, we showed that there is a large amount of symmetry or repetition in the coefficients of these polynomials.

%\subsection{}
\end{document}